\documentstyle{article}
\title{Some adjunction-theoretic properties \\ of codimension 
two nonsingular subvarieties of quadrics}
\author{Mark Andrea A.  de Cataldo}
\date{2/5/1996}

\newtheorem{tm}{Theorem}[section]
\newtheorem{lm}[tm]{Lemma}
\newtheorem{pr}[tm]{Proposition}
\newtheorem{rmk}[tm]{Remark}
\newtheorem{cor}[tm]{Corollary}

\newtheorem{?}[tm]{Question}
\newtheorem{conj}[tm]{Conjecture}

\newtheorem{fact}[tm]{Fact}

\font\tenmsb=msbm10
\font\sevenmsb=msbm7
\font\fivemsb=msbm5

\newfam\msbfam
\textfont\msbfam=\tenmsb
\scriptfont\msbfam=\sevenmsb
\scriptscriptfont\msbfam=\fivemsb
\def\Bbb#1{{\fam\msbfam #1}}

\font\teneufm=eufm10
\font\seveneufm=eufm7
\font\fiveeufm=eufm5
\newfam\eufmfam
\textfont\eufmfam=\teneufm
\scriptfont\eufmfam=\seveneufm
\scriptscriptfont\eufmfam=\fiveeufm
\def\frak#1{{\fam\eufmfam\relax#1}}

\newcommand\ci{\cite}

\newcommand\s{\sigma}

\newcommand\zed{{\Bbb Z}}
\newcommand\pn[1]{{\Bbb P}^{#1}}

\newcommand{\Q}[1]{{\cal Q}^{#1}}

\newcommand\blacksquare{{\hspace*{\fill} $\Box$}} 
\newcommand\odix[1]{ {\cal O}_{#1} }
\newcommand\odixl[2]{ {\cal O}_{#1}({#2}) }

\newcommand\nb[2]{ {\cal N}_{ {#1}, {#2} } }

\begin{document}
\maketitle

\begin{abstract}
We make precise the structure of the  first two reduction morphisms
associated with codimension two nonsingular subvarieties
 of quadrics $\Q{n}$,  $n\geq 5$.
We give a coarse classification of the same class of subvarieties 
when they are assumed
to be not of log-general-type.

\end{abstract}

\section{INTRODUCTION}
\label{intr}

Because of the Barth-Larsen Theorem and the Double Point Formula,
 low codimensional embeddings in projective space
are special in many respects.
Inspired by the study of the  special adjunction-theoretic properties
of threefolds in $\pn{5}$ contained  in 
\ci{be-sc-soams}, in this note we study the similar properties for 
codimension two  nonsingular subvarieties of quadrics $\Q{n}$, $n\geq 5$.
As it turns out, by analogy with the results of \ci{be-sc-soams},
the reduction morphisms associated with these varieties are almost 
always
isomorphisms; see Theorem \ref{secondreduction}. We give a coarse 
classification Theorem
for the varieties  for which the second reduction morphism is not 
defined, the so-called
 varieties  {\it not of log-general-type}; see Theorem 
\ref{k+2lnotspanned}, Theorem
\ref{k+2lspannednotbig}
and Theorem \ref{k+nu-2lnotnefbig}. To prove the latter one we need  
to  analyze the case 
of Del Pezzo fibrations and,  in the same way as in 
the paper \ci{boss2}, 
 the case of conic bundles on $\Q{5}$; see sections 
\ref{fanofibrations} and
\ref{quadricfibration}, respectively.

\smallskip \noindent
{\bf Notation and conventions.}
Our basic reference is [Ha].
We work over any algebraically closed field of characteristic zero.
A quadric $\Q{n}$,
here, is a nonsingular hypersurface of degree two in the  projective 
space
$\pn{n+1}$.
Little or no distinction is made between line bundles, associated 
sheaves of sections
and Cartier divisors.

\noindent
By {\em scroll} we mean a variety
$X \subseteq  \pn{N}$, for 
which $(X,\odixl{\pn{N}}{1}_{|X})  \simeq ({\Bbb P}_Y({\cal E}),
 \xi_{\cal E})$,
where $\cal E$ is a vector bundle on a nonsingular
variety $Y$. An adjunction-theoretic scroll (see \ci{be-fa-soadjtp})
is not, in general, a scroll; we denote them by 
{\em a.t. scrolls}.

\smallskip
\noindent
{\bf Acknowledgments.} This paper is an expanded
and completed version of parts of our dissertation.
It is a pleasure to thank our Ph.D. advisor
A.J. Sommese, who has suggested to us that we  study threefolds
on $\Q{5}$. We thank the C.N.R. of the Italian Government
and The University of Notre Dame for partial support.

\section{PRELIMINARY MATERIAL}

Let $\iota: X \hookrightarrow \Q{n}$ be the embedding of a degree $d$ 
 nonsingular subvariety of codimension two of 
$\Q{n}$; let $L$ denote the line bundle $\iota ^*\odixl{\Q{n}}{1}$, 
$g$ the genus of the curve $C$ obtained by intersecting
$(n-3)$ general elements of $|L|$. 
Denote by $x_i$
the Chern classes of the tangent bundle of $X$ and by
$n_i$
the ones of the normal bundle $\nb{X}{\Q{n}}$; by adjunction
$K_X=-nL +n_1$ and by the self-intersection formula
$n_2=(1/2)dL^2$.

The  following formul\ae\, which hold in the Chow ring of $X$ for 
$n\geq 5$,   are obtained using
the Double Point Formul\ae \ (see \ci{fu}) for  $\iota$. 
 
\begin{equation}
n_2=\frac{1}{2}(n^2 -n +2)L^2 -n x_1\cdot L + x_1^2 - x_2;
\label{deg2dpf}
\end{equation}
\begin{equation}
\frac{1}{6}(n^3 -3n^2+8n-12)L^3 +\frac{1}{2}(-n^2 +n -2)x_1 L^2+
n(x_1^2-x_2)L +2x_1x_2-x_1^3 -x_3=0.
\label{deg3dpf}
\end{equation}

\noindent
The following formula for surfaces $X$ on $\Q{4}$
with balanced cohomology class can be found in \ci{a-s}.

\begin{equation}
\label{deg2dpf4}
2K_X^2=\frac{1}{2}\,d^2-3d-8(g-1) +12\chi ({\cal O}_X).
\end{equation}

In the case of $n=5$,
using the above formulae we can express
$K_X\cdot L^2$, $K_X^2\cdot L$, $K_X^3$, $x_2 \cdot L$ and $x_3$ as 
functions
of $d$, $g$, $\chi ({\cal O}_X)$, $\chi ({\cal O}_S)$; for example, 
omitting the dots
from now on:

\begin{equation}
\label{KL2}
K_X L^2=2(g-1)-2d,
\end{equation}
\begin{equation}
\label{K2L}
K_X^2 L=\frac{1}{4}d^2 + \frac{3}{2}d -8(g-1) +6\chi ({\cal O}_S),
\end{equation}
\begin{equation}
\label{K3}
K_X^3=-\frac{9}{4}d^2+ \frac{27}{2}d+gd +18(g-1)-
30\chi ({\cal O}_S)-24\chi ({\cal O}_X).
\end{equation}

\bigskip

\begin{pr}
\label{s3>=0forN}
Let $X$ be a nonsingular threefold on $\Q{5}$.
Then
$$
60 \chi (\odix{S}) \geq \frac{3}{2}d^2 - 12d + (d-48)(g-1) + 24
\chi (\odix{X})
$$
and
$$
\chi (\odix{S})\leq \frac{2}{3} \frac{(g-1)^2}{d} - \frac{1}{24}d^2 + 
\frac{5}{12}d.
$$
\end{pr}

\noindent
{\em Proof.}
Denote by $s_i$ and $n_i$ the Segre and Chern classes respectively
 of the normal bundle 
$\cal N$ of $X$ in $\Q{5}$. Since $\cal N$ is generated by global 
sections,
we have $s_3\geq 0$. Since $s_3=n_1^3- 2n_1n_2$, we get
$$
0\leq (K_X+5L)^3 - 2(K_X+5L)\frac{1}{2}dL^2=
K^3 + 15K_X^2L + 75K_XL^2 + 125d - d(K_X+5L)L^2.
$$
The first inequality follows from
(\ref{K3}), (\ref{K2L}) and (\ref{KL2}).

\noindent
We use   the Generalized Hodge Index Theorem 
of \ci{boss1} (see also \ci{be-bi-so}):
$$
d(K_X^2L)\leq (K_XL^2)^2
$$
and we make explicit the left hand side using
(\ref{K2L}) and the right hand side using
(\ref{KL2}). The second inequality follows.
\blacksquare

\medskip

In what follows:

 - $((a,b,c),{ \cal O}(1))$ denotes the polarized pair given by a  
complete intersection
of type  $(a,b,c)$ in $\pn{n+1}$ and the restriction of the 
hyperplane bundle to it;

- $(X,L)$  denotes the polarized pair given by
a variety $X\subseteq \Q{n}$ and $L:= \odixl{\Q{n}}{1}_{|X}$;

- $g,$ $q$ and  $p_g$  denote the sectional genus of the embedding
line bundle, the irregularity and geometric genus of a 
surface section, respectively.

\begin{rmk}
\label{degeven}
Let $X\subseteq \Q{n}$, $n\geq 5$, be any subvariety. 
Then the degree $d$ of $X$ is even.
This follows from the fact that the cohomology class of $[X]$ equals
the class $(1/2)\, d\, [\Q{n-2}]$
in $H^4(\Q{n}, \zed)$. 
\end{rmk}

\begin{pr}
\label{classificationd<12}
{\rm (Cf. \ci{decascroll})}
Let $X\subseteq \Q{n}$, $n\geq 5$, a codimension 
two nonsingular subvariety of degree $d\leq 10$. Then the pair
$(X,L)$ is one of the ones below.

\smallskip
\noindent
\underline{\rm Type A):} $d=2$,   $((1,1,2), {\cal O}(1))$; 
$g=q=p_g=0$.

\noindent
\underline{\rm Type B):} $d=4$,   $((1,2,2),{\cal O}(1))$;  
$g=1$, $q=p_g=0$.

\noindent
\underline{\rm Type C):} $d=4$,    $n=6$, 
$(\pn{1}\times \pn{3}, {\cal  O}(1,1))$;
 $g=q=p_g=0$. 

\noindent
\underline{\rm Type D):} $d=4$,    $n=5$, 
$({\Bbb P}( {\odixl{\pn{1}}{1}}^2\oplus
\odixl{\pn{1}}{2}), \xi)$; \quad {\rm 5)}; $g=q=p_g=0$.

\noindent
\underline{\rm Type E):} $d=6$,   $((1,2,3),{\cal O}(1))$;
\quad $g=4$, $q=0$, $p_g=1$.

\smallskip
\noindent
\underline{\rm Type F):} $d=6$,  $n=5$, 
$({\Bbb P}( {\cal T}_{\pn{2}}),\xi)$, 
embedded using a general
codimension one linear system ${\frak l}\subseteq| 
\xi_{ {\cal T}_{\pn{2}} }|$;
$g=1$, $q=p_g=0$.

\smallskip
\noindent
\underline{\rm Type G):}  $d=6$
  $n=5$, $f: X \to \pn{1} \times \pn{2}=:Y$
a double cover, branched
along a divisor of type $\odixl{Y}{2,2}$,
$L\simeq   p^* \odixl{ Y }{1,1}$;  $g=2$, $q=p_g=0$.

\noindent
\underline{\rm Type H):} $d=8$,   $((1,2,4), {\cal O}(1))$; 
\quad $g=9$, $q=0$, $p_g=5$.

\smallskip
\noindent
\underline{\rm Type I):}  $d=8$, $((2,2,2), {\cal O}(1) )$; 
$g=5$, $q=0$, $p_g=1$.

\smallskip
\noindent
\underline{\rm Type L):}  $d=8$, 
$n=5$, $({\Bbb P}( E), \xi)$, $E$ a rank two vector bundle on 
$\Q{2}$ as in
{\rm \ci{io3}};  $g=4$, $q=p_g=0$.

\noindent
\underline{\rm Type M):} $d=10$,  $((1,2,5), {\cal O}(1) )$; 
$g=16,$ $q=0,$ $p_g=14$.

\noindent
\underline{\rm Type N):} $d=10$, $n=5$, $f_{|K_X+L|}: X \to \pn{1}$ 
is a fibration
with general fiber a Del Pezzo surface $F$, $K_F^2=4$, $K_X=-L+
f^*\odixl{\pn{1}}{1}$;  $g=8$, $q=0$, $p_g=2$. 

\end{pr}

We say that a nonsingular threefold $X$
on $\Q{5}$ is of Type O), if it  has degree $d=12$ and it
is a scroll over a minimal $K3$ surface. Such a threefold exists. 
See \ci{decascroll}.

\begin{pr}
\label{maintm}
{\rm (Cf. \ci{decascroll})}
The following is the complete list of \,nonsingular codimension 
two subvarieties
of quadrics $\Q{n}$, $n\geq 5$, which are scrolls.

\noindent
{\rm Type C)}, $n=6$, $d=4$,  scroll over $\pn{1}$ and over $\pn{3}$;
 
\noindent
{\rm Type D)}, $n=5$, $d=4$, scroll over $\pn{1}$; 

\noindent
{\rm Type F)}, $n=5$, $d=6$, scroll over $\pn{2}$;

 \noindent
{\rm Type \,L)}, $n=5$, $d=8$, scroll over $\Q{2}$;

\noindent
{\rm Type O)},  $n=5$, $d=12$, scroll over a minimal $K3$ surface.
\end{pr}

\begin{pr}
\label{coreasybound}
{\rm  
(Cf. \ci{decacurves}, or \ci{a-s} for the case
$d> 2k(k-1)$.)
}
Let $C\subseteq \Q{3}$ be an integral curve of degree $d$ and 
geometric genus
$g$. Assume that $C$ is contained in a surface of $\Q{3}$ of 
degree $2k$.  
Then
$$
g-1\leq \frac{d^2}{4k} + \frac{1}{2} (k-3)d.  
$$
\end{pr}

\begin{pr}
\label{boundasep} 
{\rm (Cf.  \ci{a-s}, Proposition $6.4$.)}
Let $C$ be an integral curve in $\Q{3}$, not contained in any
surface of $\Q{3}$ of degree strictly less than $2k$. Then:
$$
g-1\leq \frac{d^2}{2k} +\frac{1}{2}(k-4)d.
$$
\end{pr}

\noindent
Let $S$ be a nonsingular surface on $\Q{4}$, ${\cal N}$ 
its normal bundle, $\sigma$
 its postulation,  $C$ a nonsingular hyperplane
section of $S$, $g$ its genus, $d$ its degree.
Let $s$ be a positive integer, $V_s \in |{\cal I}_{S,\Q{4}}(s)|$
be integral and 
$\mu_l:=c_2({\cal N}(-l))=$ $(1/2)d^2 +l(l -3)d-2l(g-1)$, 
$\forall l\in \zed$.

\begin{lm}
\label{epas}
In the above situation:
\quad $
0\leq \mu_{s}\leq s^2 d.
$

\end{lm}

\noindent
{\em Proof.}
The left hand  side inequality is just Proposition
  \ref{coreasybound} above.
To prove the right hand side we first assume
$s=\s$.
Using \ci{a-s}, Lemma 6.8 we conclude (from here on the hypothesis
$d>2\s^2$ was not used there) in the case at hand.

\noindent
Now, for the general case, let
$s=\s + t$, where $t$ is a non-negative integer. Then, as it is easily
checked,
$\mu_s=\mu_{\s}+ \s td + t(\s+t-3)d -2t(g-1)$. We conclude by
what proved for $\mu_{\s}$ and by  the obvious
$g\geq 0$.
\blacksquare

\begin{rmk}
\label{barthlarsen}
{\rm
Let $X$ be a nonsingular codimension two subvariety of 
$\Q{n}$. 
As a consequence of the Barth-Larsen Theorem (see
\ci{ba}),  we have that: 
if $n\geq 6$, then
 the fundamental group $\pi_1(X)$ is trivial; 
if $n\geq 7$,  then $Pic(X)\simeq \zed$, generated by the hyperplane 
bundle,
so that $X$ does not carry any nontrivial morphisms.
}
\end{rmk}

\medskip
The following fact  is well known when $\Q{n}$ is replaced by
$\pn{n}$, see  \ci{boss1} for example. The case
of $\Q{4}$ is proved in \ci{a-s}, Lemma 6.1. The general case
can be proved in the same way.  See
\ci{bounded}, where we prove a more general statement.
We used this ``lifting" criterion  as a tool to prove the finiteness of 
the number of families of nonsingular
threefolds on $\Q{5}$ not of general type; see Proposition
\ref{sigh} below. 
 
\begin{pr}
\label{roth}
{\rm (Cf. \ci{bounded})}
Let $X$ be an integral subscheme of degree $d$ and codimension two on 
$\Q{n}$,
 $n\geq 4$. Assume that for the general hyperplane section
$Y$ of $X$ we have
$
h^0({\cal I}_{Y,\Q{n-1}}(\s))$ $\not= 0,
$
for some positive integer $\s$ such that $d>2{\s}^2$.
Then\,
$
h^0({\cal I}_{X,\Q{n}}({\s}))\not= 0.
$
\end{pr}

\begin{pr}
\label{sigh}
{\rm (Cf. \ci{bounded})}
Let $n=4,$ $5$ or $n\geq 7$.
There are only finitely many components of the Hilbert scheme
of $\Q{n}$ corresponding to nonsingular $(n-2)$-folds not of general 
type.
\end{pr}

\section{THE STRUCTURE OF THE REDUC\-TION MOR\-PH\-IS\-MS}
\label{dpfand divisors}

In this section we give, by a systematic use of the 
double point formul\ae, 
 a precise description of the reduction morphisms 
associated with  codimension two subvarieties of quadrics $\Q{n}$, 
$n\geq 5$. 
We apply these formul\ae \ also
to the case of divisorial contractions of extremal rays on 
threefolds on $\Q{5}$.
For the language and results of Adjunction Theory, which we 
are going to use freely
for the rest of this note,
 we refer the reader to
\ci{be-fa-soadjtp} and to \ci{be-so}.

\smallskip
Let  $\nu :=n-2$.

\begin{lm}
\label{numericaldpf}
Let $X$ be a codimension two nonsingular subvariety of $\Q{n}$, 
$n\geq 5$.

\noindent
Let $D$ be a divisor on $X$ with $(D,\odixl{D}{D})\simeq (\pn{\nu-1},
\odixl{\pn{\nu -1}}{-1})$
and $(K_X+(\nu-1)L)_{|D}\simeq \odix{D}$;
then $n=5,$ $6$ and $d=10$.

\noindent
Let $n=5$. Then we have the following list of possible degrees
according to whether $X$ contains a divisor of the given form 
 $(D,\odixl{D}{D})$ with   $(K_X+(\nu-2)L)_{|D}\simeq \odix{D}$:

\smallskip
\noindent
{\rm (\ref{numericaldpf}.1)} 
if  $(D,\odixl{D}{D}) \simeq ( \pn{2}, \odixl{ \pn{2}  }{-2})$, then
$d=20$;

\noindent
{\rm (\ref{numericaldpf}.2)}  
if $(D, \odixl{D}{D}) \simeq  (\pn{2}, \odixl{\pn{2}}{-1})$, then
$d=14$;

\noindent
{\rm (\ref{numericaldpf}.3)}
 if  $(D ,\odixl{D}{D}) \simeq ( \tilde{ {\Bbb F}}_2, G)$, where 
$2G=K_D$, then
$d=14$;

\noindent
{\rm (\ref{numericaldpf}.4)} 
$(D,\odixl{D}{D}) \simeq (  {\Bbb F}_0,  G)$, where $2G=K_D$, 
then $d=14$;

\noindent
{\rm (\ref{numericaldpf}.5)} the case in which $D$  has two 
components as in
 {\rm \ci{be-fa-soadjtp}, Theorem 0.2.1, case b5)}, 
 cannot occur;

\noindent
{\rm (\ref{numericaldpf}.6)}  the case 
$(D ,\odixl{D}{D}) \simeq ( {\Bbb F}_1, -E-f)$ cannot occur.

\noindent
{\rm (\ref{numericaldpf}.7)} the cases in which $D$ is as in either \,
{\rm a)}, or {\rm  b)} 
of \,   {\rm  \ci{be-so}}, {\rm  Theorem 2.3} cannot occur.

 \smallskip
\noindent
Let $n=6$. Assume $X$ contains a surface $\cal S$ such that
${\cal S}\simeq \pn{2}$, $L_{|{\cal S}} \simeq \odixl{\pn{2}}{1}$ and 
such that
the normal bundle ${\cal N}_{{\cal S},X} \simeq 
{\cal T}^*_{ \pn{2} } (1)$.
Then $d=14$. 

\end{lm}

\noindent
{\em Proof.} For $n=5$ the proof is the same  as the one 
of \ci{be-sc-soams}, Proposition 1.1, using
(\ref{deg2dpf}) in the place of (0.8) of the quoted paper.
For $n=6$ we compute all the relevant Chern classes
 by using (\ref{deg2dpf}),  the Euler sequence
for ${\cal S}\simeq \pn{2}$ and the exact sequence
$$
0\to {\cal T}_{\cal S} \to {{\cal T}_X}_{|{\cal S}} 
\to {\cal N}_{{\cal S},X} \to 0.
$$
\blacksquare

\begin{tm}
\label{secondreduction}
{\rm (\bf{Structure of the reduction morphisms})}
Let $X$ be a nonsingular codimension two subvariety of $\Q{n}$, 
$n\geq 5$.

\noindent
Assume that $(X,L)$ admits a first reduction $(X',L')$.
Then the first reduction morphism is an isomorphism: 
$(X,L)\simeq (X',L')$.

\noindent
Assume that $(X,L)$ admits, in addition, a second reduction 
$(X'',L'')$. We have:

\noindent
if \ $n=5$ and $d\not= 14,$ $20$, then $(X,L)=(X',L')$ and the 
second
reduction map $\varphi: X'\to X''$
 is the blowing up on a nonsingular $X''$ of a disjoint union
of nonsingular integral curves; 

\noindent
if  $n=6$ and $d\not= 14$,
then $(X,L)=(X',L')$ and the second
reduction map $\varphi: X'\to X''$
 is the blowing up on a nonsingular $X''$ of a disjoint union
of nonsingular integral curves. 
If in addition $d\not=
16,$ $22$,  then the second reduction morphism is an isomorphism:
$(X,L)\simeq$ $(X',L')\simeq$ $(X'',L'')$;

\noindent
if $n\geq 7$,  then $(X,L)\simeq $$(X',L')\simeq$$(X'',L'')$.
\end{tm}

\noindent
{\em Proof.}

Once $K_X+(n-1)L $ is nef and big, i.e. out of the lists of Theorem
\ref{k+2lnotspanned} and Theorem \ref{k+2lspannednotbig}, it 
fails to be ample
only if the first reduction is not an isomorphism; in turn, 
that happens
if and only if $X$ contains  some exceptional divisors of the 
first kind. By
Proposition \ref{numericaldpf} this happens only if $d=10$. 
By Proposition
\ref{classificationd<12}
the  type is either
M) or N); neither of them contains an exceptional divisor of 
the first kind. It follows that
 if the first reduction exists,
then $(X,L)\simeq (X',L')$.

\noindent
The statements concerning the second reduction morphism can
be proved as follows. For $n=5,$ we use  Theorem
0.2.1 of  \ci{be-fa-soadjtp}
coupled with Proposition \ref{numericaldpf}. 

\noindent
For
$n=6$ we use Theorem 0.2.2 of 
\ci{be-fa-soadjtp} and then
we take a general hyperplane  section and reduce to the case
$n=5$, with the difference that now  case b2) of Theorem 0.2.1
of \ci{be-fa-soadjtp} does not occur.
 The case of  the blowing up of curves
yields $d=16,\,22$, as we now show.
Since $X\simeq X'$ we cut (\ref{deg2dpf})
with $F\simeq \pn{2}$, a general fiber of the blowing up.
Define $a$ to be the positive integer such that
$L_{|F}\simeq \odixl{\pn{2}}{a}$. Since
$\nb{F}{X}\simeq$ $\odix{\pn{2}} \oplus$ $\odixl{\pn{2}}{-1}$ and
${K_X}_{|F} \simeq$ $\odixl{\pn{2}}{-2}$ we get
$$
(16-d/2)a^2=12a -4.
$$
Since $a>0$ we see that $d\leq 30$. The only integer solutions to the
relation above are
$(d,a)=$ $(16,1)$ and $(22,2)$. This concludes the proof
for $n=6$.

\noindent
 Finally, for $n\geq 7$ we use Remark \ref{barthlarsen}.
\blacksquare

\medskip

We now describe  Mori contractions for threefolds on $\Q{5}$.

\begin{lm}
\label{divisorialmoricontractionsd=}
Let $X$ be a nonsingular threefold in $\Q{5}$. Let $D$ be an 
integral divisor
on $X$. We have:

\smallskip
\noindent
{\rm (\ref{divisorialmoricontractionsd=}.1)} 
if $(D,\odixl{D}{D}) \simeq (\pn{2}, \odixl{\pn{2}}{-1})$, then either
$d=10$ and $L_{|D}\simeq \odixl{\pn{2}}{1}$, or $d=14$ and
$L_{|D}\simeq \odixl{\pn{2}}{2}$;

\noindent
{\rm (\ref{divisorialmoricontractionsd=}.2)} 
if $(D,\odixl{D}{D}) \simeq (\pn{2}, \odixl{\pn{2}}{-2})$, then either
$d=8$ and $L_{|D}\simeq \odixl{\pn{2}}{1}$, or $d=16$ and
$L_{|D}\simeq \odixl{\pn{2}}{2}$;

\noindent
{\rm (\ref{divisorialmoricontractionsd=}.3)}
 if $(D,\odixl{D}{D}) \simeq ( {\Bbb F}_0,G)$, then
$d\leq 20$;

\noindent
{\rm (\ref{divisorialmoricontractionsd=}.4)} 
if $(D,\odixl{D}{D}) \simeq ( \tilde{\Bbb F}_2,G)$, then $d=14$ and
$L_D=-G$.

\end{lm}

\noindent
{\em Proof.}  The proof is the same  as the one of \ci{be-sc-soams},
 Proposition 1.1, using
(\ref{deg2dpf}) in the place of (0.8) of the quoted paper.
\blacksquare

\begin{pr}
\label{moricontrd>=22}
{\rm (\bf{Structure of Mori contractions})}
Let $X$ be a nonsingular threefold embedded in $\Q{5}$ with
 $d\geq 22$ and $K_X$
not nef. Let $\rho: X\to Y$ be the contraction of any extremal 
ray  on
$X$. Then $Y$ is nonsingular and either $\rho$ is birational 
and the blowing up of
an integral nonsingular  curve on $Y$ or $\rho$ is a conic 
bundle in the sense
of Mori Theory. In particular, if $d\gg 0$, then only the 
former case can occur.
\end{pr}

\noindent
{\em Proof.}  The proof is the same  as the one of \ci{be-sc-soams}, 
Corollary 1.2, using
(\ref{deg2dpf}) in the place of (0.8) of the quoted paper.
As for the last statement, if $\dim Y \leq 2$, then $X$ is 
not of general type
and we apply Proposition \ref{sigh}
\blacksquare

\medskip
The following conjecture is due to Beltrametti, Schneider 
and Sommese in the case 
of $3$-folds on $\pn{5}$. It seems a fairly
natural question in view of Proposition \ref{moricontrd>=22}.

\begin{conj}
\label{conjbescso} There is an integer $d_0$ such that 
every threefold
on $\Q{5}$ of degree $d\geq d_0$ is a minimal model.

\end{conj}

\section{VARIETIES NOT OF LOG-GENERAL-TYPE}
\label{secadjointbundles12}
\noindent

In this section we give a coarse classification of varieties as in 
the title.
We still make free use of the language of Adjunction Theory. 

\smallskip
Let $(X,L)$ be a degree $d$, $\nu$-dimensional nonsingular 
subvariety of $\Q{n}$ endowed
with its embedding line bundle $L$.
 The ``Types"  we shall consider correspond to the ones of 
Propositions 
\ref{classificationd<12} and \ref{maintm}.

\smallskip
We start by observing that $K_X+(\dim X-1)L$ is spanned by 
its global sections
(spanned for short)
except for three varieties.

\begin{tm}
\label{k+2lnotspanned}
Let $(X,L)$ be as above.  Then $K_X+(\nu -1)L$ is 
spanned  unless
$(X,L)$ is one of the three pairs {\rm A)}, {\rm C)} or {\rm D)}.
In particular,
$d\leq 4$.
\end{tm}

\noindent{\em Proof.} By the list on \ci{be-so} page 381,
  and by the fact
that there are no codimension two linear subspaces on 
${\Q{n}},$ $\forall n\geq 5$,
 we need to analyze the a.t. scroll  over a curve case only.
By flatness an a.t. scroll over a curve is a scroll.
The result follows from Theorem \ref{maintm}.
\blacksquare

\medskip

Now we classify those pairs for which $K_X+(\nu -1)L$ is spanned, 
but for which 
$\kappa (K_X+(\nu -1)L)<\nu$.

\begin{tm}
\label{k+2lspannednotbig}
Let $(X,L)$ be as above. Assume that $K_X+(\nu -1)L$
is spanned, i.e. $(X,L)$ is not as in
{\rm Theorem   \ref{k+2lnotspanned}}, but that it is not big. Then 
$(X,L)$ is one of the following pairs:

\noindent
{\rm (\ref{k+2lspannednotbig}.1) (Del Pezzo variety):} {\rm Type B)}; 
{\rm Type F)};

\noindent
{\rm (\ref{k+2lspannednotbig}.2) (Quadric Bundle over a curve):}
{\rm Type G)};

\noindent
{\rm (\ref{k+2lspannednotbig}.3) (A.t. scroll over a surface):}
{\rm Type L)}; {\rm Type O)}.

\noindent
In particular, $d\leq 12$.

\end{tm}

\noindent
{\em Proof.}
Let $K_X+(\nu -1)L$ be as in the Theorem, then
by \ci{be-so} page 381 $(X,L)$ is either a Del Pezzo variety, 
a quadric bundle or an
a.t.  scroll
over a surface.

\noindent
Let us assume that $(X,L)$ is a Del Pezzo variety.
By slicing with $(\dim X-2)$ general hyperplanes we get
a surface in $\Q{4}$ with $K_S=-L_{|S}$. Since $S$ is Del Pezzo
we get $\chi (\odix{S})=g(L)=1$. 
We plug these values in 
(\ref{deg2dpf4}) and get:
$$
d^2-10d+24=0.
$$
It follows that either $d=4$ or $d=6$. The conclusion follows from 
Proposition
\ref{classificationd<12}.

\noindent
Let us assume that $(X,L)$ is a quadric bundle. 
Let $F\simeq \Q{n-3}$ be a general fiber of the  quadric fibration. 
Dotting
(\ref{deg2dpf}) with $F$ we get $d=6$. We conclude using Proposition
\ref{classificationd<12}.

\noindent
Let us assume that $(X,L)$ is an a.t. scroll over a surface.
By \ci{be-so-book}, Proposition 14.1.3
 $(X,L)$ is an ordinary scroll with $\kappa (K_X + (n-1)L)=2$.
We conclude by comparing with 
Proposition \ref{maintm}.
\blacksquare

\medskip
Now we deal with the line bundle $K_X+(\nu -2)L$.
 First we exclude the presence of some special pairs.

\begin{lm}
\label{noveronese}
Let $(X,L)$ be as above, then $(X,L)$
cannot be isomorphic to any of the three pairs
$(\pn{4}, \odixl{\pn{4}}{2})$, $(\pn{3},\odixl{\pn{3}}{3})$ and
$(\Q{3},\odixl{\Q{3}}{2})$. Moreover,
there are no Veronese bundles $(X',L')$ associated with a pair
$(X,L)$ on $\Q{5}$.
\end{lm}

\noindent
{\em Proof.}
By contradiction assume that $(X,L)\simeq (\pn{4}, \odixl{\pn{4}}{2})$. 
We intesect two general members of $|L|$ and get
 a  nonsingular surface section $(S,L_{|S})$ which is embedded
in $\Q{4}$ with $d=16$, $g=1$ and $\chi(\odix{S})=1$. This  contradicts
(\ref{deg2dpf4}). We exclude the case in which
$(X,L)\simeq (\Q{3},\odixl{\Q{3}}{2})$ in a similar way.

\noindent
The possibility $(X,L)\simeq (\pn{3},\odixl{\pn{3}}{3})$
is ruled out by Remark \ref{degeven}.

\noindent
Let us assume that $(X,L)$ is a pair for which $(X',L')$ exists
and is a Veronese bundle with
associated morphism $p:X \to Y$; in particular $n=5$.
By Theorem \ref{secondreduction} $(X,L)\simeq (X',L')$. Dotting
(\ref{deg2dpf}) with a general fiber
$F$ we get $d=10$.  Since for some ample line bundle $\cal L$ on $Y$
$2K_X+3L=p^*{\cal L}$,  we have 
the following relation on   a general surface section $S$  of $X$:
$$
L_{|S}=-2K_S + {\cal L}_{|S},
$$
which ``squared" gives $d=10 \equiv 0$ $mod(4)$, a contradiction.
\blacksquare

\begin{tm}
\label{k+nu-2lnotnefbig}
 Assume that we are not on the lists of {\em Theorems 
\ref{k+2lnotspanned}}
and {\em   \ref{k+2lspannednotbig}} so that $(X,L)\simeq (X',L')$. 
If $K_X+(\nu -2)L$
is not nef and big then $(X,L)$ is one of the following pairs:

\noindent
{\rm (\ref{k+nu-2lnotnefbig}.1) (Mukai variety):} {\rm Type E)}; {\rm
Type I)};

\noindent
{\rm (\ref{k+nu-2lnotnefbig}.2) (Del Pezzo fibration over a curve):}
either {\rm Type N)}, $d=10$ or as in 
{\rm (\ref{3dimdelpezzofibr}.2)}, $d=12$; 

\noindent
{\rm (\ref{k+nu-2lnotnefbig}.3) (Quadric bundle over  a surface):}
$n=5,\,6$,
 a flat quadric bundle over a nonsingular surface: if $n=6$, then $d=12$
and if $n=5$, then  either $d\leq 18$  or $d=44$.

\noindent
{\rm (\ref{k+nu-2lnotnefbig}.4) (A.t. scroll over a threefold):}
$n=6$, the scroll map  is not flat and $d$ is either $14$ or $20$.

\end{tm}
\noindent
{\em Proof.}
Let $K_X+(\nu -1)L$ be as in the Theorem, then
by \ci{be-so} page 381-2    and   Lemma \ref{noveronese},
 $(X,L)$ is either 
 a Mukai variety, a Del Pezzo fibration over a curve, a 
quadric bundle over a surface
or an a.t. scroll of dimension $\nu \geq 4$ over a normal threefold.

\noindent
Let us assume that $(X,L)$ is a Mukai variety.
By slicing to a surface section $S$ we find that $K_S=\odix{S}$, 
and since
$X$ is simply connected it follows that $\pi_1(S)$ is trivial as well;
$S$ is thus a $K3$ surface. Using (\ref{deg2dpf4}) we get,
using $\chi (\odix{S})=2$, $2(g-1)=d$, that either $d=6$ or $d=8$; 
accordingly
$g=4,\,5$, respectively. The conclusion, in this case, follows from
Proposition \ref{classificationd<12}.

\smallskip
\noindent
We deal with the case of Del Pezzo fibrations over a curve in Lemma 
\ref{delpezzofibrationandkodaira} and Proposition
\ref{3dimdelpezzofibr}

\smallskip
\noindent
 We now deal with quadric bundles over surfaces.
Again, $n=5,\, 6$, by Remark \ref{barthlarsen}.

\noindent
Let $n=5$ and assume, by contradiction, that there is a divisorial 
fiber $F$ of the
quadric bundle map $p:X \to Y$.
Then $F$ is as in \ci{be-so} Theorem 2.3. This contradicts case 
(\ref{numericaldpf}.7) of Lemma \ref{numericaldpf}. It follows 
that all the fibers 
of $p$ are 
equidimensional. By  Theorem
\ref{conbundle} it follows that $p$ is a quadric fibration in
the sense
of section \ref{quadricfibration}. The statement follows form
Proposition
\ref{conicbundleonsurface} and Remark \ref{maple}.

\noindent 
Let $n=6$.  $(X,L)$ is a quadric bundle over a surface, $p:X \to Y$, 
so is
its general  hyperplane section. By what proved for the case $n=5$ 
the base surface $Y$ is nonsingular and  by Corollary
\ref{liftflat}
we deduce that $p$ is flat. If we cut (\ref{deg2dpf}) with a general 
fiber
of $p$ we get $d=12$.

\noindent
Case (\ref{k+nu-2lnotnefbig}.3) follows.

\smallskip
\noindent
Finally case  (\ref{k+nu-2lnotnefbig}.4) follows from 
Proposition \ref{maintm} which ensures us of the absence, on $\Q{6}$, 
of
 adjunction theoretic scrolls over threefolds 
for which the map $p$ is flat: for if
$p$ were flat then $Y$ would be nonsingular by
\ci{mats} Theorem $23.7$ and then $X$ would be a projective bundle, 
a contradiction.  
 If one of these
scrolls occurs, since $p$ is not flat and -$K_X$ is
$p$-ample, Lemma \ref{conbundle} and 
\ci{mats}, Theorem $23.1$ ensures there must be  a fiber $F$ such 
that either
$F$ contains a divisor or, by \ci{be-so-book}, 14.1.4,  
$F$ is a surface $\cal S$ as in 
Proposition \ref{numericaldpf}. In the latter case we get $d=14$. 
In the former,
by slicing with a general hyperplane section, we get a threefold
$\tilde{X}$ together with the morphism
$\tilde{p}:=p_{|{\tilde X}}:{\tilde X} \to Y$, where $Y$ is the
 base of the scroll. $\tilde{p}$ 
is the second reduction morphism
for  $(\tilde{X}, L_{|{\tilde X}})$, so that the result follows 
by looking at the divisorial
fibers of $\tilde{p}$ and Lemma \ref{numericaldpf}.
\blacksquare.

\section{FIBRATIONS  OVER CUR\-VES WI\-TH GE\-NE\-RAL FI\-BER 
A DEL PEZ\-ZO MA\-NI\-FOLD}
\label{fanofibrations}

In this section we study codimension two nonsingular subvarieties
of $\Q{n}$, $n\geq 5$, which admit a   morphism
$f:X\to Y$, with connected fibers,
 onto a nonsingular curve  $Y$, such that 
the line bundle $K_X + (n-4) L$ is trivial on the general fiber. 
The general fiber will thus be a 
nonsingular (adjunction-theoretic) Del Pezzo
variety of the appropriate dimension $n-3$.
 By  Remark \ref{barthlarsen}  we have
$n=5,$ $6$.

 The following
lemma ensures that these fibrations coincide with the Del Pezzo 
fibrations
over curves  of Adjunction Theory.

\begin{lm}
\label{delpezzofibrationandkodaira}
Let $X$ be a fibration as above. Then $K_X+(n-1)L$ is ample and
$\kappa (K_X+(n-2)L)=\kappa (S)=1$.
\end{lm}

\noindent
{\em Proof.}
\noindent
Without loss of generality we may assume that   $n=5$, for
otherwise we cut with a general hyperplane section to the
three dimensional case and it is easy to show that  if the 
statements we
 want to prove hold
for the threefold hyperplane section  of $X$, then they also 
hold for $X$.
 
\noindent
The generic fiber of $f$ is a nonsingular Del Pezzo surface
$F$. Since $K_X+L$ is  trivial on the fibers
we define
$$
\Delta:=L^2\cdot F=L_{|F}^2=K_F^2.
$$
Cut (\ref{deg2dpf}) with $F$, using the facts that
${K_X}_{|F}=K_F$ and that $x_2\cdot F=12- \Delta$. We get
$$
\Delta=\frac{24}{16-d}.
$$
Since $F$ is a Del Pezzo surface and $L$ is very ample, we get
$3\leq \Delta\leq 9$. 
Since $\Delta$ is an integer we have only the following possibilities:
\begin{equation}
\label{deltaandd}
(\Delta, d)= (3,8),\ (4,10), \ (6,12).
\end{equation} 
Using the above invariants, and the lists of Adjunction Theory, 
it is easy to show that  $K_X+(n-1)L$  is  ample
and that
 $\kappa(K_X+(n-2)L)=0,$ $1$.
By Theorem \ref{k+nu-2lnotnefbig}
 the case $K_X=-(n-2)L$ cannot occur, 
since these manifolds do not carry any
nontrivial fibration. It follows that 
$K_X+2L$ is ample, $\kappa(K_X+L)=1$ and, by adjunction,
$\kappa(S)=1$.
\blacksquare 

\medskip

We need the following facts.

\begin{fact}
\label{relvan}
{\rm 
Let $f:X\to Y$ be as above. By relative vanishing we have
 $h^i(\odix{X})=h^i(\odix{Y})$, $\forall i$.
}
\end{fact}

\begin{fact}
\label{socrelle}
 $g(Y)=q(S)$, $2g-2-d=(p_g(S) + q(S) -1)\Delta$; 
 moreover the elliptic fibration
$S\to Y$ has no multiple fibers.
\end{fact} 

\noindent
The assertion on $g(Y)=q(S)$ follows from Lefschetz 
Theorem on hyperplane sections,
$q(S)=h^1(\odix{X})$, 
and from Fact \ref{relvan}; the other assertion  follows from 
\ci{socrelle}, 0.5.1.

\begin{fact}
\label{notinp4}
  $ S\not \subseteq  \pn{4} $.
\end{fact}

\noindent
 To prove this, assume that
$S\subseteq {\pn{4}}$.
We use jointly  the double point formula for surfaces on 
$\pn{4}$, see \ci{ha}, page 434,
and (\ref{deg2dpf4}) to compute the values of
$g$ and $\chi (\odix(S))$ to conclude that, $d=8,$ $10$ 
would yield
noninteger values, a contradiction,  and
 that if $d=12$ then $g=25$, and $\chi (\odix{S})=13$;
this system of invariants is inconsistent by  Fact
\ref{socrelle}. This proves the assertion.

\begin{pr}
\label{3dimdelpezzofibr} Let $X\subseteq \Q{n}$, $n\geq 5$,
be a nonsingular, 
codimension two subvariety
which  admits a
 fibration
$f:X\to Y$ in Del Pezzo manifolds onto   a nonsingular curve $Y$; 
in particular
$(K_X+(\dim X -2)L)_{|F}\simeq \odix{F}$, $F$ a  general fiber.

\noindent
Then $Y\simeq \pn{1}$ and either $(X,L)$ is of {\rm Type N)}
or only the following systems of invariants is possible:

\noindent
{\rm (\ref{3dimdelpezzofibr}.1)} 
$n=5,$ $d=12$, $K_F^2=6$, $g=10$, $p_g(S)=2$, $q(S)=0$,  
$h^i(\odix{X})=0,\ \forall i>0$.
\end{pr}

\noindent
{\em Proof.} By the  proof of Lemma
\ref{delpezzofibrationandkodaira} and by the knowledge
of degree $d=8,\,10$ varieties stemming from
Proposition \ref{classificationd<12}, we only need to rule 
out the case $d=8$
and  make precise   the invariants in the case
$d=12$. Moreover, by the same lemma,  $\kappa (S)=1$.

\noindent 
First let $n=5$.

\noindent
Now we determine the invariants in the case $d=12$.

\noindent
 We apply  formula (\ref{deg2dpf4}) 
in the case $d=12$. We get
\begin{equation}
\label{chi}
2(g-1)-3\chi (\odix{S})=9.
\end{equation}
 By Fact \ref{notinp4} and by Proposition
\ref{roth} we are in the position to apply the Castelnuovo bound for
curves on $\pn{4}$, which   gives
$g\leq 13$.

\noindent
 (\ref{chi}) implies that $\chi (\odix{S})$ is not a non-negative 
integer, unless
$(g,\chi (\odix{S}))=$  $(7,1)$, $(10,3)$, $(13,5)$.
We can rule out the cases:  $d=12$ and $(g,\chi(\odix{S}))=$$(7,1)$,  
$(13,5)$
using Fact \ref{socrelle}   which gives $g-7=3(p_g + q -1)$; 
this last equality  
together with the given values
of $\chi(\odix{S})$ and $g$ gives a non-integer value for $q$, a 
contradiction. 
It follows that
 if $d=12$,
then  $(g,\chi (\odix{S}))=$ $(10,3)$.
 To compute the values of $p_g$ and $q$ 
we use again Fact \ref{socrelle} which gives the number
$p_g+q$. Since we know $\chi(\odix{S})$ we get the values of 
$p_g$ and $q$.

\noindent
 Since $g=q$ we see that $Y\simeq \pn{1}$. The assertions on 
$h^i(\odix{X})$
follow from Fact \ref{relvan}.

\noindent
The proposition is thus proved for  $n=5$.

\noindent
Let $n=6$, the only remaining case, by Barth-Larsen Theorem. 
By slicing with a general hyperplane we get a threefold
with a fibration onto a curve whose general fiber is a Del Pezzo
 manifold so 
that the above analysis applies. The only difference is that the 
case  $d=10$
does not occur by Proposition \ref{classificationd<12}.

\noindent
Now we  prove that also the case $d=12$ does not occur. 

\noindent
The general fiber of $f$ is a Del Pezzo threefold with
$K_F=-2L_{|F}$ and $L_{|F}^3=6$. By explicit classification,
 see \ci{fujbook}, page 72,
either
$F\simeq \pn{1}\times \pn{1} \times \pn{1}$ or $F\simeq 
{\Bbb P}({\cal T}_{\pn{2}})$.
In both cases formula (\ref{deg3dpf}) dotted with $F$ 
gives $x_3\cdot F=x_3(F)=24$.
But in the former  case $x_3(F)=8$, in the latter $x_3(F)=6$.
\blacksquare

\subsection{MORE UPPER BOUNDS}
This section is not needed for Theorem \ref{k+nu-2lnotnefbig}.

\medskip
 We now give an upper bound for the degree of codimension two, 
nonsingular
subvarieties of $\Q{n}$, $n\geq 5$, which admit a morphism onto
a curve such that the  a 
general fiber is a Fano variety. By Barth-Larsen Theorem
we need to worry only about the cases $n=5,$ $6$.

\begin{pr}
\label{boundforfanofibrations}
Let $X\subseteq \Q{n}$ a nonsingular subvariety of codimension 
two and degree $d$
 which
admits a morphism onto a curve such that the general fiber is a 
Fano variety.

\noindent
If $n=5$, then $d\leq 20$.

\noindent
If $n=6$, then $d\leq 30$.

\end{pr}

\noindent
{\em Proof.}
Let $n=5$ and ${\cal L}:=L_{|F}$. Assume that $d\geq 22$.

\noindent
 We cut (\ref{deg2dpf}) with a fiber, $F$, and obtain, on $F$:
$$
(11-d/2){\cal L}^2 + 5K_F{\cal L} + K_F^2 -c_2(F)=0.
$$
Since  $c_2(F)=12- K_F^2$,  we get:
\begin{equation}
\label{=ondelpezzo}
(d/2-11){\cal L}^2  + 2K_F^2 -12 + 5K_F {\cal L}=0.
\end{equation}
Now we use $K_F^2\leq 9$ to get
\begin{equation}
\label{fibrdelpezzosurfaces}
(d/2 -11){\cal L}^2\leq  6 + 5K_F{\cal L}.
\end{equation}
Since $K_F{\cal L}\leq -1$, we see that either $d=22$, or $d=24$ and
${\cal L}^2=- K_F{\cal L}=1$. In the latter case $F\simeq \pn{2}$ and 
 the Hodge Index Theorem on the surface $F$ says that $K_{\pn{2}}^2=1$,
 a contradiction. In the former case
we use (\ref{=ondelpezzo}):
$$
2K_F^2-12 + 5K_F{\cal L}=0,
$$
which gives a contradiction for each value
$K_F^2=1,\ldots,$ $9$. It follows that $d\leq 20$.

\noindent
The  proof of the statement for  $n=6$ is analogous to the proof
of Proposition \ref{fanod<22}, where we use (\ref{deg2dpf}) with
$n=6$ cut with the cycle $K_X\cdot F$.
\blacksquare

\medskip

 In the same spirit we give an upper bound on the degree of Fano
threefolds on $\Q{5}$.

\begin{pr}
\label{fanod<22}
Let $X\subseteq \Q{5}$ be a nonsingular  Fano threefold. 
Then $d\leq 20$.
\end{pr}

\noindent
{\em Proof.} (Cf. \ci{be-sc-soams}, Corollary 1.2.)
We cut (\ref{deg2dpf}) with $K_X$ and get, using the fact that
$x_1x_2=24 \chi (\odix{X})=24$:
$$
(11-d/2)L^2K_X + 5LK_X^2 + K_X^3 + 24=0.
$$
Let
$$
\lambda:= LK_X^2, \qquad 2\mu:=- L^2K_X=-2g+2+2d;
$$
clearly $\lambda$ and $\mu$ are positive integers and the above becomes:
\begin{equation}
\label{fanolm}
(d-22)\mu + 5\mu \lambda + 24 = -K_X^3.
\end{equation}
By the Generalized Hodge Index Theorem, see \ci{boss1}, we get
$(-K_X^3)(-K_XL^2)\leq$ $ (K_X^2L)^2$, or
\begin{equation}
\label{hodgeonfano}
(-K_X)^3(2\mu)\leq \lambda^2.
\end{equation}
By combining 
(\ref{fanolm}) and (\ref{hodgeonfano}) we get
\begin{equation}
\label{l2m}
\lambda^2 - 10\mu \lambda - [2(d-22)\mu^2 + 48\mu]\geq 0.
\end{equation}
If we solve the above in $\lambda$ we get either $\lambda <0$, a 
contradiction,
or $\lambda > 10 \mu$. This implies, in turn, that $\lambda \geq 11$.
Since, by the classification of Fano threefolds, $-K_X^3\leq 64$, 
(\ref{fanolm}) becomes
$$
(d-22)\mu +55 + 24\leq 64,
$$
a contradiction for $d\geq 22$.
\blacksquare

\section{QUADRIC FIBRATIONS}
\label{quadricfibration}

In this section the term ``quadric bundle" is to be intended in 
the sense of Adjunction Theory.
The term ``quadric fibration" is introduced below.

By  {\em quadric fibration}
we mean a nonsingular projective variety $X\subseteq {\Bbb P}$, of
dimension $x$, together with a fibration
$p:X\to Y$ onto a ({\em a fortiori}) 
nonsingular variety $Y$ of positive dimension $y$, {\em all}
of which fibers are quadrics, not necessarily
integral, of the appropriate dimension $(x-y)$.
One has non integral fibers
only if the relative dimension is one.

The case $\dim Y=0$ is trivial. 
In virtue of Remark \ref{barthlarsen} we have:

\begin{fact}
{\em There are no codimension two quadric fibrations on $\Q{n}$, 
for $n\geq 7$ and,
for $n=6$, any such is simply connected.}
\end{fact}

\medskip
We restrict ourselves to the case of $n\geq 5$.

We begin by fixing some notation and establishing some simple facts.

\noindent
Let  $L$ denote the restriction to $X$ of the hyperplane bundle. 
The sheaf
${\cal E} : =p_* L$ is  locally free on $Y$ of rank $(x-y +2)$. 
It is easy to check that $\cal E$ is generated by  its global sections.
The surjection
$p^*p_*L\to L$
defines an embedding:
$X \hookrightarrow {\Bbb P}(\cal E)$, where $L={\xi_{\cal E}}_{|X}$
and $X$ is defined by a nonzero section of the line bundle
$ 2\xi - \pi^* {\cal M}$, for some ${\cal M}\in Pic (Y)$, where 
$\pi:{\Bbb P}({\cal E})\to Y$ is the bundle projection. 

The following gives a sufficient condition for a general hyperplane 
section
of $X$ to be a quadric fibration over $Y$. 
 It is a well known ``counting dimensions" argument.

\begin{lm}
\label{countingdimensions}
Let $X\to Y$ be a quadric fibration as above. Assume
$2y< x+2$. Then a general hyperplane section $X'$ of $X$ is a quadric 
fibration
over $Y$ via $p_{|X'}:X'\to Y$.  
\end{lm}

\noindent
{\em Proof.} Since $\cal E$ is generated by global sections and, by 
assumption
$rank({\cal E})>y$, a general section of it does not vanish on $Y$. 
Such a section
will define, {\em for every} $y\in Y$, a hyperplane $\Lambda_y$ of
the corresponding fiber
$\pi^{-1}(y)\subseteq {\Bbb P}({\cal E})$. In the case in which
the quadrics $p^{-1}(y)$ were integral $\forall y\in Y$,
we would be done. This is, in general,  not true. 
However, the singular quadrics of the fibration are parameterized
by a proper closed subset $D$ of $Y$ with $\dim D \leq (y-1)$.
The hyperplanes of $\Bbb P$ which contain the reduced 
part, $\Sigma\simeq \pn{x-y}$, of one of the components of one non 
integral 
quadric of the fibration
form a linear space  of dimension $(\dim {\Bbb P}-x+y-1)$
contained in ${\Bbb P}^{\vee}$. The space of these bad hyperplanes 
is of
dimension at most $(\dim D + \dim {\Bbb P} -x +y-1)\leq$ $\dim 
{\Bbb P}-x+2y-2<$
$\dim {\Bbb P}^{\vee}$. It follows that  the general section
of $\cal E$ gives a hyperplane section of $X$ which cuts {\em every}
quadric of the fibration  in a quadric of dimension one less.
\blacksquare  

\medskip

\begin{pr}
\label{quadricfibrationoncurve}
There are no quadric fibrations over curves on $\Q{6}$.
The only quadric fibrations over curves on $\Q{5}$ are of 
{\rm Type G)}. 
If there is a quadric fibration over a surface on $\Q{6}$, 
hen it has degree $d=12$.
 \end{pr}

\noindent
{\em Proof.} As to quadric fibrations over curves, we
 cut (\ref{deg2dpf}) with a nonsingular fiber
$F\simeq \Q{n-3}$, we  get $d=6$. We conclude by comparing with
Proposition \ref{classificationd<12}.

\noindent
As to quadric fibrations over a surface we cut 
(\ref{deg2dpf}) with  a nonsingular fiber $F\simeq \Q{n-4}$ and 
get $d=12$.
\blacksquare

\medskip

The following proposition and remark describe the situation for 
threefolds
 quadric bundles over surfaces.

\begin{pr}
\label{conicbundleonsurface}
Let $X\subseteq \Q{5}$ be a threefold quadric fibration
{\rm (conic bundle)} over a surface $Y$. Then either 
$d\leq 98$ or $X$ is contained in a hypersurface $V\in 
|\odixl{\Q{5}}{3}|$ and $d\leq 276$.
\end{pr}

\noindent
{\em Proof.} We denote the Chern classes of $X$ and $Y$ by 
$x_i$ and $b_i$, 
respectively.
We omit the symbol ``$p^*$" for ease of notation.
 We follow closely the paper \ci{boss2}. First we introduce the 
following
entities and we report from \ci{boss2}, for the reader's convenience,
the relations among them which are essential to the computations below
(one warning:
 some of the equalities are only numerical equalities):

$\cal M$ was defined at the beginning of the section;

$D\in |2e_1 -3{\cal M}|$, it is called the {\em discriminant divisor}; 
its points correspond to  the singular fibers of $p$;

$2R\subseteq Y$ the {\em branching divisor} associated with a general
 hyperplane section, $S$, of
$X$, which, in view of Lemma \ref{countingdimensions}, is a cyclic 
double
cover of $Y$;
 
$e_1= 3R-D$;

${\cal M}=2R-D$;

$x_1=L +b_1 - R;$

$x_2=L^2 +L\cdot( b_1 -2R +D) + ( -2R^2 - R\cdot b_1+ D\cdot R  + 
b_2 + e_2) ;$

$x_3= 2b_2 -D^2 + Db_1;$

$L\cdot W\cdot W'=2W\cdot W'$, for every pair of divisors $W$ and 
$W'$ on $Y$;

$L^2\cdot W= (4R-D)\cdot W$;

$e_2=\frac{1}{2}(12R^2 +D^2 -7DR -d)$.

Now we plug the above values of $x_1$ and $x_2$ for $x_1$ and $x_2$ in
(\ref{deg2dpf}):
\begin{equation}
\label{deg2cb}
(6-\frac{d}{2})L^2 -4L b_1+5L R+ b_1^2 -b_1R  
-LD + 3R^2-D R
-b_2 -e_2=0.
\end{equation}
Next we equate the expression above for $x_3$ to the one
of
(\ref{deg3dpf}), using again the above expressions for $x_1$ and $x_2$:
\begin{equation}
\label{deg3cb}
-(2d+10)b_1R + 2dR^2
+ (\frac{d}{2}+4)Db_1 + D^2 -10b_2 
+ 2b_1^2 -(\frac{d}{2}+5)DR -d(\frac{d}{2} -13)=0.
\end{equation}
Now we set 
$$
x:= b_1^2 \quad {\rm and} \quad y:=DR,
$$
 we cut $(\ref{deg2cb})$ 
with $R$, $-b_1$, $D$ and $L$, respectively, so that we obtain
four linear equations to which we add (\ref{deg3cb}), after having 
substituted
in $x$ and $y$. The result is the following linear system of equations:

\begin{equation}
\label{systcb}
Mv^t=c^t,
\end{equation}
where

$
\hspace{1in}  M:= \left( 
\begin{array}{rrrrr}
-8 & 34-2d & 0 & 0 & 0                  \\
2d-34 & 0 & -\frac{d}{2}+8 & 0 & 0      \\
0 & 0 & -8 & \frac{d}{2}-8 & 0           \\
-18 & 14 & +4 & 0 & -2                  \\
-2d-10 & 2d & \frac{d}{2} +4 & 1 & -10
\end{array}
\right),
$

\smallskip
\noindent
 $v:=(\ b_1R,\ R^2,\  Db_1,\  D^2,\  b_2\ )$

\noindent
and 

\noindent
$c:=(\  (8-\frac{d}{2})y,\  -8x,\  (2d-34)y,\  2x+
4y + d(\frac{d}{2}-7),\ 
-2x + (\frac{d}{2} +5)y + d(\frac{d}{2}-13)\ )$.

\medskip
\noindent
Since
$P:=-\frac{1}{2}det M=3d^3-27d^2 -1520d +18976>0$, $\forall d>0$, 
we can solve the above
system
(\ref{systcb}) and obtain the unique solution:

\medskip
\noindent
$b_1R =  -\frac{1}{2}[(-128 d^2 +4480d-39168)x+(2d^3 -111d^2+2020d 
-12096)y+$ $(
2d^5 -120d^4 +2678d^3 -26304d^2+95744d)]/P,$

\medskip
\noindent
$R^2 = \frac{1}{4}[(-1024d+18432)x +(3d^3 -8d^2 -2112d+23552)y +
(16d^4 -688d^3 +9728d^2 -45056d)]/P, $    

\medskip
\noindent
$b_1D = -2[(-152d^2 +4440d-32128)x + (2d^3 -113d^2 +2099d -12852)y 
+ (
2d^5 -122d^4 +2766d^3 -27574d^2 +101728d)]/P$, 

\medskip
\noindent
$D^2 =-4[(-1216d +16064)x + (-3d^3 +46d^2 +893 d -13736)y +
(16d^4 -720 d^3 +10608d^2 -50864d)]/P, $       

\medskip
\noindent
$b_2 =\frac{1}{4}[(12d^3+20d^2 -3648d+13952)x+(d^3 -30d^2 +152d+960)y 
+
(d^5 -27d^4 +274d^3 -4448d^2 +46016d)]/P $.

\medskip
\noindent
Since $\cal E$ is generated by global sections and $D$ is effective
we see that
$e_2\geq 0$,
 $e_1D\geq 0$.
Also, \ci{boss2}, Lemma 
2.9 gives
$y=DR\geq 0$.
 We can make explicit 
$e_2$ and $e_1$ 
by the formul\ae \ 
given at the beginning 
of this proof and deduce:
\begin{eqnarray*}
\label{ineqcb}
DR &= & y  \geq  0,          \\
e_2 \cdot P & = & (896d-4480)x - (\frac{19}{2}d^2 -366d+3616)y-     \\
& & (\frac{19}{2}d^4                          
-\frac{843}{2}d^3 + 5864d^2 -24656d)\geq 0,     \\
 e_1D \cdot P &= & -(4864d-64256)x- (3d^3 -103d^2 + 988d -1984)y+     \\
& & (64d^4 -2880d^3 + 42432d^2 -203456d)\geq 0.
\end{eqnarray*} 
These three inequalities define a region of the plane $(x,y)$.
It is straightforward to check that the two lines
$e_2=0$ and $e_1D=0$ have slopes $a$ and $b$ whose sign
does not change with $d$ if $d\geq 20$. One can check easily that
$a>0$ and $b<0$. The intersection of the first line above with the 
$x$-axis
is
$$
(x_1,0)_{e_2}=(\frac{(19/2)d^4 -
(843/)2d^3+5864d^2 - 24656d}{896d-4480},0);
$$
the intersection of the second line with the $x$-axis is
$$
(x_2,0)_{e_1D}=(\frac{64d^4 -2880d^3 + 
42432d^2 -203456d}{4864d-64256},0).
$$
One can check, that, since 
$d\geq 20$, 
$x_1<x_2$.
The region we are interested in is a triangle with vertices
$(x_1,0)_{e_2}$,
$(x_2,0)_{e_1D}$ 
and 
$(x_3,y_3)_{(e_2=0)\cap (e_1D=0)}.$

\noindent
Now we compute the genus of a general curve section, $C$, of $X$.  
By adjunction
$x_1\cdot L^2=2d+2-2g$, so that by what above:

\begin{eqnarray*}
g-1 & =& \frac{d}{2} -2b_1R+ \frac{Db_1}{2} + 2R^2 -\frac{DR}{2}   \\
 &= & -2b1R +  \frac{Db_1}{2} +2R^2 - \frac{y}{2} +\frac{d}{2}    \\
&=&[\,(24d^2-472d+2176)x+( (23/2)d^2 -375d +3044)y +   \\
& & ((23/2)d^4 -(891/2)d^3 +5374 d^2 -19024d)\,]\,/\,P.
\end{eqnarray*}
Again it is not difficult to check that the absolute value of the slope
of the above line is bigger than 
$|b|$. 
It follows easily that the maximum possible value for
$g-1$ 
in our region is achieved at
$(x_2,0)_{e_1D}$,
 while the minimum is at
$(x_1,0)_{e_2}$.
We thus get
\begin{equation}
\label{superbound}
\frac{19d^3-187d^2+416d}{224d-1120} \leq g-1
\leq \frac{4d^3-77d^2+321d}{38d-502}.
\end{equation}
Assume that $C$ is not contained in any surface of $\Q{3}$ of 
degree strictly less
than $2\cdot 11$. Then by (\ref{boundasep}) and by the left hand side
inequality of (\ref{superbound}), we get
$$
\frac{19d^3-187d^2+416d}{224d-1120} \leq  \frac{d^2}{22}+\frac{7}{2}d,
$$
which, remembering  that $d$ is even and
that we are assuming $d\geq 20$, implies
$ d\leq 98$.

\noindent
Assume that 
$C$ is contained in a surface of degree $2k$, with $k=10,9,\ldots,
3$. By Corollary (\ref{coreasybound})
we infer:
$$
\frac{19d^3-187d^2+416d}{224d-1120} \leq  \frac{d^2}{4k}+\frac{k-3}{2}d,
$$
which implies, as above, that
for $k=10, 9, \ldots, 3$,
$d\leq 64,$ $58$ $54,$ $48,$ $44,$ $40,$ $40$  and $276$, respectively.

\noindent
Finally, assume that $C$ is contained in a surface of degree four or
 two.
Using the right hand side inequality of (\ref{superbound})
and Lemma \ref{epas} we get
$d\leq  42$ and $d\leq 16 $, respectively. Actually in the last case 
we get a contradiction,
since we are assuming $d\geq 20$.

\noindent
Finally if $C$ is in a surface of degree six, then $X$ is in a 
hypersurface
of degree six in $\Q{5}$, provided, $d>18$ (cf. Proposition \ref{roth}).
\blacksquare

\begin{rmk}
\label{maple}
{\rm We have checked with a Maple routine which are the possible 
degrees of
a threefold on $\Q{5}$ which is
a quadric fibration over a surface.  For $d\geq 20$
 we have imposed the following restrictions on the
triples $(d,x,y)$:

\noindent
1) $20 \leq d \leq 276$;

\noindent
2) for every fixed $d$ as above $(x,y)$ must belong to the triangle
of the proof of Proposition \ref{conicbundleonsurface};

\noindent
3) $b_1R,$ $  R^2,$ $ b_1D,$ $D^2,$  $b_2$, $g-1$, $\chi ({\cal O}(Y)$  
and
$\chi ({\cal O} (S))$ must be integers; 

\noindent
4) $(g-1)$ must satisfy inequality (\ref{superbound}) and the bound
of Theorem 2.3 in \ci{gross};

\noindent
5) $\chi ({\cal O}(S))$ must satisfy the two inequalities
 of
Proposition \ref{s3>=0forN};

\noindent
6) various inequalities stemming from the Hodge Index Theorem on $Y$
as, for example,  $(K_YR)^2\geq K_Y^2 R^2$;

\noindent
7) if $d> 98$ then $g-1 \leq (1/12)d^2$, see Proposition 
\ref{coreasybound};

\noindent
The result is that the only possible degree, for $d\geq 20$ is $d=44$.

\noindent
By taking double covers of the four scrolls of \ci{ottp5},
we see that there are flat conic bundles over surfaces for 
$d=6,\, 12,\, 14,\, 18$.
 We do not know 
whether the case $d=44$ occurs.
}  
\end{rmk}

\subsection{DIGRESSION}
\label{digressione}

In the course of the proof of Theorem \ref{k+nu-2lnotnefbig}
we used the fact, due to Besana
\ci{besana}, that the base of an adjunction theoretic quadric 
bundle over a surface is nonsingular. 
The following is a result with a similar flavor. It is probably
well known.

\begin{lm}
\label{conbundle}
Let $X$ a nonsingular projective variety of dimension $n$, $p:X\to Y$
a morphism onto a normal projective variety $Y$ of dimension $n-1$ 
such that all fibers
have the same dimension, 
the general scheme theoretic fiber over a closed point  is 
isomorphic to a conic and
$-K_X$ is $p$-ample.
Then all the scheme-theoretic fibers are isomorphic to conics, 
$p$ is flat
and $Y$ is nonsingular.
\end{lm}

\noindent
{\em Proof.}
The proof is the same as the one of \ci{mori} Lemma 3.25.
The only necessary changes are the following: a) replace the 
line bundle
$H$ of \ci{mori}, by a pull-back $p^*A$ of any ample line bundle 
$A$ on $Y$ and
use
Kleiman  criterion of ampleness to obtain the result analogue to 
the last assertion of 
\ci{mori} Lemma;
b) replace\ci{mori} Lemma 3.12 by \ci{ando} Lemma 1.5.
\blacksquare

\begin{cor}
\label{liftflat}
Let $X$ be a nonsingular projective variety together with a morphism
$p:X\to Y$, where $Y$ is a normal variety of  dimension $m$.
Let $D_i$, $i=1,\ldots, n-m-1$ be divisors on $X$ such that they 
intersect
 transversally; denote by  $X'$  their intersection. Assume that 
$p_{|X'}:X'\to Y$ satisfies the hypothesis of {\rm  Lemma 
\ref{conbundle}}.
Then $p$ is flat and $Y$ is nonsingular.
\end{cor}

\noindent
{\em Proof.} By the lemma, $p_{|X'}$ is flat. We can ``lift" 
this flatness
to $p$ by virtue of \ci{mats}, Corollary to  Theorem 22.5. As above
the flatness of $p_{|X'}$ (or of $p$) implies the nonsingularity of
$Y$.
\blacksquare

\begin{cor}
\label{conics}
Let $X$ a nonsingular projective variety of dimension $n$, $p:X\to Y$
a morphism onto a normal projective  variety $Y$ of dimension $n-1$ 
such that all fibers
have the same dimension.
If the general fiber of $p$ is actually embeddable as conics with
 respect to 
an embedding of $X$, then all scheme theoretic
fibers are actually embedded conics, $p$ is flat, $Y$ is nonsingular 
and
$-K_X$ is $p$-ample.
\end{cor}

\noindent
{\em Proof.} We argue as in the proof of the lemma with the 
simplifications
due to the fact that a flat  deformation of a conic in projective space
 is still a conic.
The assertion on $-K_X$ follows by observing that,  if $L$ denotes 
the line bundle with which we embed $X$, $K_X+L$ is a pull-back
from $Y$.
\blacksquare 

\begin{rmk}
{\rm The assumption $-K_X$ is $p$-ample is essentialin the lemma, 
as the blow up of a $\pn{1}$
bundle over a curve at two distinct points on a fiber shows.
Moreover, the above Lemma does not follow directly from
\ci{mori} or  \ci{ando}, since there are conic bundles which 
structural morphism
is not a  Mori contraction. 
Finally,
the above theorem is certainly 
false if one has $\dim X= \dim Y$. It is a purely local question:
consider the quotient of $\Bbb A^2$ by the involution
$(x,y)\to (-x,-y)$.
}
\end{rmk}

\bigskip
1991 {\em Mathematics Subject Classification}. 14C05, 14E05, 14E25, 
14E30, 14E35, 14J10,
14J30, 14J35, 14J40, 14J45, 14M07.

{\em Key words and phrases}. Adjunction Theory, classification, 
codimension two,
conic bundles, low codimension, non log-general-type, quadric, 
reduction, 
special variety.

\bigskip
AUTHOR'S ADDRESS

\medskip
\noindent
Mark Andrea de Cataldo,
Department of Mathematics,
Washington University in St. Louis,
Campus Box 1146,
St. Louis, Missouri 63130-4899.

\noindent
e-mail: mde@math.wustl.edu

\end{document}